\documentclass{amsart}  %used for final

\usepackage{amsmath}%,cancel,fullpage}
\usepackage{tikz}
\usepackage[all]{xy}
\usetikzlibrary{decorations.markings}
\usepackage{color}
\usepackage{amsthm}
\usepackage{amsfonts}
\usepackage{amssymb}
\usepackage{setspace}
\usepackage{youngtab}
\usepackage{ytableau}
\usepackage{empheq}
\usepackage{comment}
\usepackage{cancel} % for gigantic cancellation marks
\usepackage[bookmarks,colorlinks,breaklinks]{hyperref}  % PDF hyperlinks, with coloured links
\hypersetup{linkcolor=blue,citecolor=blue,filecolor=blue,urlcolor=blue} % all blue links
\theoremstyle{plain}

\theoremstyle{plain}
\newtheorem{theorem}{Theorem}[section]

\newtheorem{lemma}[theorem]{Lemma}
\newtheorem{corollary}[theorem]{Corollary}

\newtheorem{definition}[theorem]{Definition}
\newtheorem{example}[theorem]{Example}

\theoremstyle{remark}
\newtheorem{remark}[theorem]{Remark}

\numberwithin{equation}{section}
\newcommand{\bN}{\mathbb{N}}

\newcommand{\coverc}{\lessdot _{c}} % less than cover relation left composition poset
\newcommand{\coverr}{\lessdot _{r}} % less than cover relation right composition poset
\newcommand{\coverq}{\lessdot _{q}} % less than cover quasisymmetric composition poset
\newcommand{\coverqt}{< _{\widetilde{q}}} % less than cover def composition poset

\newcommand{\lessc}{< _{c}} % less than relation left composition poset
\newcommand{\lessr}{< _{r}} % less than relation right composition poset
\newcommand{\lessq}{< _{q}} % less than relation quasisymmetric composition poset
\newcommand{\Lc}{\mathcal{L}_{c}} % left Composition poset

\newcommand{\Rc}{\mathcal{R}_c} % right Composition poset
\newcommand{\Qc}{\mathcal{Q}_c} % quasisymmetric Composition poset
\newcommand{\Qct}{\widetilde{\mathcal{Q}}_c} % def quasisymmetric Composition poset

\newcommand{\addt}{\mathfrak{t}} % box-adding operator
\newcommand{\addu}{\mathfrak{u}} % jdt operator

\newcommand{\wD}{\widetilde{D}} % down operator for dual filtered graph
\newcommand{\Ut}{U_t} % up operator for left composition poset
\newcommand{\Id}{{Id}} % identity operator 
\newcommand{\down}{\mathfrak{d}}
\newcommand{\downrow}{\mathfrak{d}}
\newcommand{\addurow}{\mathfrak{u}}

\newcommand{\suchthat}{\;|\;}

\newlength\cellsize 
\savebox2{%
\begin{picture}(15,15)
\put(0,0){\line(1,0){15}}
\put(0,0){\line(0,1){15}}
\put(15,0){\line(0,1){15}}
\put(0,15){\line(1,0){15}}
\end{picture}}
\newcommand\cellify[1]{\def\thearg{#1}\def\nothing{}%
\ifx\thearg\nothing
\vrule width0pt height\cellsize depth0pt\else
\hbox to 0pt{\usebox2\hss}\fi%
\vbox to 15\unitlength{
\vss
\hbox to 15\unitlength{\hss$#1$\hss}
\vss}}
\newcommand\tableau[1]{\vtop{\let\\=\cr
\setlength\baselineskip{-16000pt}
\setlength\lineskiplimit{16000pt}
\setlength\lineskip{0pt}
\halign{&\cellify{##}\cr#1\crcr}}}
\savebox3{%
\begin{picture}(15,15)
\put(0,0){\line(1,0){15}}
\put(0,0){\line(0,1){15}}
\put(15,0){\line(0,1){15}}
\put(0,15){\line(1,0){15}}
\end{picture}}
\newcommand\expath[1]{%
\hbox to 0pt{\usebox3\hss}%
\vbox to 15\unitlength{
\vss
\hbox to 15\unitlength{\hss$#1$\hss}
\vss}}
\newcommand\bas[1]{\omit \vbox to \cellsize{ \vss \hbox to \cellsize{\hss$#1$\hss} \vss}}

\begin{document}
\title[Dual graphs from generalized Schur functions]{Dual graphs from noncommutative and quasisymmetric Schur functions}
\author{S. van Willigenburg}
\address{Department of Mathematics, University of British Columbia, Vancouver, BC V6T 1Z2, Canada}
\email{\href{mailto:steph@math.ubc.ca}{steph@math.ubc.ca}}
\thanks{The author was supported in part by the National Sciences and Engineering Research Council of Canada.}
\begin{abstract}
By establishing relations between operators on compositions, we show that the posets of compositions arising from the right and left Pieri rules for noncommutative Schur functions can each be endowed with both the structure of dual graded graphs and dual filtered graphs when paired with the poset of compositions arising from the Pieri rules for quasisymmetric Schur functions and its deformation.
\end{abstract}

\maketitle

\section{Introduction}\label{sec:intro} Differential posets \cite{stanley-diff} and dual graded graphs \cite{fomin-dual1,fomin-dual2} were first developed in order to better understand the Robinson-Schensted-Knuth algorithm. However, since then they have developed into a research area in their own right, for example \cite{miller,stanley-zanello}, including rank variants \cite{stanley-var} and signed analogues \cite{lam-sign}. They also arise in the study of representations of towers of algebras \cite{bergeron-lam-li, gaetz}, have been generalized to planar binary trees \cite{nuts}, Kac-Moody algebras \cite{berg-saliola-serrano,lam-shimozono}, quantized versions \cite{lam}, and most recently to related to K-theory via dual filtered graphs \cite{patrias-pyl}. The classic example of dual graded graphs is Young's lattice paired with itself. Young's lattice appears in a variety of areas, such as being used to describe the Pieri rules for Schur functions. From this perspective, natural generalizations of Young's lattice exist arising from Pieri rules for the Schur function generalizations known as quasisymmetric Schur functions, and noncommutative Schur functions. In particular, quasisymmetric Schur functions \cite{QS} are a nonsymmetric generalization of Schur functions that form a basis for the increasingly ubiquitous Hopf algebra of quasisymmetric functions, for example \cite{DKLT, gessel, hersh-hsiao}. Their Pieri rules \cite[Theorem 6.3]{QS} give rise to the generalization of Young's lattice known as the quasisymmetric composition poset. Dual to this Hopf algebra is the Hopf algebra of noncommutative symmetric functions  \cite{GKLLRT}, whose basis dual to that of quasisymmetric Schur functions is the basis of noncommutative Schur functions \cite{BLvW}, a noncommutative analogue of Schur functions. Due to noncommutativity, two sets of Pieri rules arise, one arising from multiplication on the right \cite[Theorem 9.3]{tewari} and one from multiplication on the left \cite[Corollary 3.8]{BLvW}. These two sets of Pieri rules give rise to two generalizations of Young's lattice known as the right composition poset and the left composition poset. Therefore the question arises: Are these  posets dual graded and dual filtered graphs? In this paper we answer this question in the affirmative.

More precisely, this paper is structured as follows. In Section~\ref{sec:comps} we review necessary notions on compositions in order to define operators on them. These operators are used to define three partially ordered sets in Subsection~\ref{subsec:posets}, $\Rc$ and $\Lc$ that arise in the right and left Pieri rules for noncommutative Schur functions, and $\Qc$ that arises in the Pieri rules for quasisymmetric Schur functions. We then establish useful relations satisfied by these operators in Subsections~\ref{subsec:ud} and \ref{subsec:td}. In Section~\ref{sec:dualgraphs} we show that $\Rc$ and $\Qc$, plus $\Lc$ and $\Qc$, are each a pair of dual graded graphs in Theorems~\ref{the:RcQc} and \ref{the:LcQc}. We define a strong filtered graph $\Qct$ on the set of compositions using the operators arising in the Pieri rules for quasisymmetric Schur functions in Definition~\ref{def:Qct}, and establish that $\Rc$ and $\Qct$, plus $\Lc$ and $\Qct$, are each a pair of dual filtered graphs in Theorems~\ref{the:RcQct} and \ref{the:LcQct}.

%Compositions

\section{Compositions and operators}\label{sec:comps} A finite list of integers $\alpha = (\alpha _1, \ldots , \alpha _\ell)$ is called a \emph{weak composition} if $\alpha _1, \ldots , \alpha _\ell$ are nonnegative, and is called a \emph{composition} if $\alpha _1, \ldots , \alpha _\ell$ are positive. Note that every weak composition has an underlying composition, obtained by removing all 0s. Given $\alpha = (\alpha _1, \ldots , \alpha _\ell)$ we call the $\alpha _i$ the \emph{parts} of $\alpha$,  and the sum of the parts of $\alpha$ the \emph{size} of $\alpha$.

%Operators

%\section{Operators on compositions}\label{sec:ops} 

Now we will recall four families of operators, each of which are indexed by positive integers, and have already contributed to the theory of quasisymmetric and noncommutative Schur functions. Although originally defined on compositions, we will define them in the natural way on weak compositions to simplify our proofs. Our first operator is the box removing operator $\down$, which first appeared in the Pieri rules for quasisymmetric Schur functions \cite{QS}. Our second operator is the appending operator $a$. Together these give our third operator, the jeu de taquin or jdt operator $\addu$. This operator arises in jeu de taquin slides on   semistandard reverse composition tableaux and in   the right Pieri rules for noncommutative Schur functions \cite{tewari}. Our fourth   operator is the box adding operator $\addt$, which plays the same role in the left Pieri rules for noncommutative Schur functions \cite{BLvW} as $\addu$ does in the   right Pieri rules. Each of these operators is defined on weak compositions for every  integer $i\geq 0$ and we note that
$$\down _0 = a_0 = \addu _0 = \addt _ 0 = \Id$$namely the identity map, which fixes the weak composition it is acting on. We now define the remaining operators for $i\geq 1$, after establishing some set notation. Let $\bN$ be the set of positive integers. Anytime we refer to a set $I\subset \bN$, we implicitly assume that $I$ is finite. Also, if we are given such a set $I$, then $I-1$ is the set obtained by subtracting $1$ from all the elements in $I$ and removing any $0$s that might arise in so doing.
\begin{example}
If $I=\{1,2,4\}$, then $I-1=\{1,3\}$.
\end{example}
By $[i]$ where $i\geq 1$, we mean the set $\{ 1,2, \ldots , i\}$.
We furthermore define $[0]$ to be the empty set.
We will denote the maximum element of a set $A$ by $\max(A)$. If $A$ is the empty set, by convention we have that $\max(A)=0$.

The first \emph{box removing operator} on weak compositions, $\down _i$ for $i\geq 1$, is defined as follows. Let $\alpha$ be a weak composition. Then
$$\down _i (\alpha) = \alpha '$$where $\alpha '$ is the weak composition obtained by subtracting 1 from the rightmost part equalling $i$ in $\alpha$. If there is no such part then we define $\down _i(\alpha) = 0$.

\begin{example}\label{ex:down}
Let $\alpha=(2,1,3)$. Then $\down_1(\alpha)=(2,0,3)$,  $\down_2(\alpha)=(1,1,3)$, $\down_3(\alpha)=(2,1,2)$ and $\down_4(\alpha)=0$. In fact, $\down_i(\alpha)=0$ for all $i\geq4$.
\end{example}

Given a finite set $I= \{i_1<\cdots <i_k\}$ of positive integers, we define
\begin{align*}
\downrow_{I}=\down_{i_1}\down_{i_2}\cdots \down_{i_k}.
\end{align*}
For convenience, we define $\downrow_{\emptyset}= \down _0$. The empty product of box removing operators is also defined to be $\down _0$. 

\begin{example}\label{ex:downset}
\begin{align*}
\downrow_{[3]}((3,1,4,2,1))&=\down_1\down_2\down_3((3,1,4,2,1))\\&=\down_1\down_2((2,1,4,2,1))\\&= \down_1((2,1,4,1,1))\\&= (2,1,4,1,0)
\end{align*}
\end{example}

The second \emph{appending operator} on weak compositions, $a_i$ for $i\geq 1$, is defined as follows. Let $\alpha = (\alpha _1, \ldots , \alpha _{\ell})$ be a weak composition. Then
$$a _i (\alpha) = (\alpha _1, \ldots , \alpha _{\ell}, i)$$namely, the weak composition obtained by appending a part $i$ to the end of $\alpha$. To simplify proofs later, we will abuse notation and also think of $a_0$ as adding 0 to the end of $\alpha$ that we will eventually remove.

\begin{example}\label{ex:append} Let $\alpha = (2,1,3)$. Then $a_2 (\alpha)= (2,1,3,2)$. However, $a_2 \down _4 (\alpha) = 0$  since $\down _4 (\alpha) = 0$ by Example~\ref{ex:down}.
\end{example}

The third \emph{jeu de taquin} or \emph{jdt operator} on weak compositions, $\addu _i$ for $i\geq 1$, is defined as follows. Considering the box removing and appending operators,
$$\addu _i = a_i \down _{[i-1]}.$$
 
 \begin{example}\label{ex:jdt} Let us compute
 $$\addu_4 ((3,1,4,2,1))=a_4\downrow_{[3]}((3,1,4,2,1)).$$By Example~\ref{ex:downset}
 $\downrow_{[3]}((3,1,4,2,1))=(2,1,4,1,0)$, and hence $\addu_4(3,1,4,2,1)=(2,1,4,1,0,4)$.
\end{example}

For any set of finite positive integers $I=\{i_1<\cdots <i_k\},$ we define
$$\addurow_I=\addurow_{i_k}\cdots\addurow_{i_1}.$$
For convenience, we define $\addurow_{\emptyset}= \addu _0$. The empty product of jdt operators is also defined to be $\addu _0$.
Note further that the order of indices in $\downrow_I$ is the reverse of that in $\addurow_I$. 

Lastly,  the fourth \emph{box adding operator} on weak compositions, $\addt _i$ for $i\geq 1$, is defined as follows. Let $\alpha = (\alpha _1, \ldots , \alpha _{\ell})$ be a weak composition. Then
$$\addt _1 (\alpha) =  (1, \alpha _1, \ldots , \alpha _{\ell})$$and for $i\geq 2$
$$\addt _i (\alpha) =  (\alpha _1, \ldots , \alpha _j + 1, \ldots ,\alpha _{\ell})$$where $\alpha _j$ is the leftmost part equalling $i-1$ in $\alpha$. If there is no such part, then we define $\addt _i (\alpha) = 0$.

\begin{example}\label{ex:boxadd}
Let $\alpha = (3,1,4,2,1)$. Then $\addt _1(\alpha)= (1,3,1,4,2,1)$, $\addt _2(\alpha)=(3,2,4,2,1)$, $\addt _3(\alpha)=(3,1,4,3,1)$,  $\addt _4(\alpha)=(4,1,4,2,1)$, $\addt _5(\alpha)=(3,1,5,2,1)$ and $\addt _i (\alpha)=0$ for all $i\geq 6$. 
\end{example}

%%%right and left and quasi posets
\subsection{Composition posets}\label{subsec:posets} With our operators we will now define three partial orders   on compositions noting that \emph{if any parts of size 0 arise during computation, then they are ignored.} The adjectives right and left in the first two are not only used to distinguish between the posets, but also to refer to their roles in the right and left Pieri rules for noncommutative Schur functions in \cite[Theorem 9.3]{tewari} and \cite[Corollary 3.8]{BLvW} respectively, and whose notation we follow now.

\begin{definition}\label{def:RcLc} The \emph{right composition poset}, denoted by $\Rc$, is the poset consisting of all compositions with cover relation $\coverr$ such that for compositions  $\alpha, \beta$
$$\beta \coverr \alpha \mbox{ if and only if } \alpha = \addu _i (\beta)$$for some $i\geq1$. Meanwhile the \emph{left composition poset}, denoted by $\Lc$, is the poset consisting of all compositions with cover relation $\coverc$ such that for compositions  $\alpha, \beta$
$$\beta \coverc \alpha \mbox{ if and only if } \alpha = \addt _i (\beta)$$for some $i\geq1$.
\end{definition}

The order relation $\lessr$ in $\Rc$ (respectively, $\lessc$ in $\Lc$) is obtained by taking the transitive closure of the cover relation $\coverr$ (respectively, $\coverc$).

\begin{example}\label{ex:RcLc} Let $\beta = (3,1,4,2,1)$, $\alpha^R = (2,1,4,1,4)$ and $\alpha^L = (4,1,4,2,1)$. Then $\beta \coverr \alpha^R = \addu _4(\beta)$ and $\beta \coverc \alpha^L = \addt _4(\beta)$ by Examples~\ref{ex:jdt} and \ref{ex:boxadd}, respectively.
\end{example}

Our third poset, meanwhile, stems from the Pieri rules for quasisymmetric Schur functions \cite[Theorem 6.3]{QS}, hence its name.

\begin{definition}\label{def:Qc} The \emph{quasisymmetric composition poset}, denoted by $\Qc$, is the poset consisting of all compositions with cover relation $\coverq$ such that for compositions  $\alpha, \beta$
$$\beta \coverq \alpha \mbox{ if and only if } \down _i(\alpha) = \beta$$for some $i\geq1$. \end{definition}

Again, the order relation $\lessq$ in $\Qc$ is obtained by taking the transitive closure of the cover relation $\coverq$.

\begin{example}\label{ex:Qc} Let $\beta = (4,1,3,2,1)$  and $\alpha = (4,1,4,2,1)$. Then $\down _4 (\alpha) = \beta \coverq \alpha$.
\end{example}

%%%u and d relations

\subsection{Relations satisfied by operators of type $\addu$ and $\down$}\label{subsec:ud} We will now give a variety of lemmas regarding the jdt operators and box removing operators, which will be useful in proving our main theorems later. Hence this subsection can be safely skipped for now and referred to later. In all the proofs we assume that $\alpha$ is a weak composition.

\begin{lemma}\label{lem:d and a same}
For $i\geq 0$ we have that $a_i=\down_{i+1}a_{i+1}$.
\end{lemma}
\begin{proof}
Let $\alpha=(\alpha_1,\ldots,\alpha_\ell)$. Then $a_{i+1}(\alpha)=(\alpha_1,\ldots,\alpha_\ell,i+1)$. This implies by definition that $\down_{i+1}a_{i+1}(\alpha)=(\alpha_1,\ldots,\alpha_\ell,i)=a_i(\alpha)$. 
\end{proof}

As a corollary we obtain the following relationship between any two appending operators.
\begin{corollary}\label{lem:sequence version d and i same}
For positive integers $i$ and $j$ satisfying $i\geq j$, we have that
\begin{align*}
\down_{j}\down_{j+1}\cdots \down_{i-1}\down_{i}a_i=a_{j-1}.
\end{align*}
\end{corollary}

\begin{lemma}\label{lem:commutativity of a and d distinct}
Let $i\neq j$ be positive integers. Then
\begin{align*}
\down_ia_j=a_j\down_i.
\end{align*}
\end{lemma}
\begin{proof}
Let $\alpha=(\alpha_1,\ldots,\alpha_\ell)$.
Let $\beta=a_j(\alpha)=(\alpha_1,\ldots,\alpha_\ell,j)$.  If $\alpha$ does not have a part equalling $i$, then neither does $\beta$, as $i\neq j$. Thus in this case we have that $\down_ia_j(\alpha)=\down _i (\beta) = 0 =a_j\down_i(\alpha)$. Now, assume that $\alpha_r$ is the rightmost part equalling $i$ in $\alpha$. Then $a_j\down_i(\alpha)=(\alpha_1,\ldots,\alpha_{r-1},\alpha_r-1,\ldots,\alpha_\ell,j)$. Since $i\neq j$, we are guaranteed that $\down_i(\beta)=(\alpha_1,\ldots, \alpha_{r-1},\alpha_r-1,\ldots,\alpha_\ell,j)$. Thus we have that $\down_ia_j(\alpha)=a_j\down_i(\alpha)$ in this case as well, and we are done.
\end{proof}

The proofs of the next three lemmas consist of case analyses that are similar in style to the proof of Lemma~\ref{lem:commutativity of a and d distinct}, however, they are more technical and hence we omit them for brevity.

\begin{lemma}\label{lem:non local commutativity for d}
Let $i$ and $j$ be distinct positive integers such that $|i-j|\geq 2$. Then 
\begin{align*}
\down_i\down_j=\down_j\down_i.
\end{align*}
\end{lemma}

\begin{lemma}\label{lem:local commutativity version 1 for d}
Let $i\geq 1$. Then $\down_i^2\down_{i+1}=\down_i\down_{i+1}\down_i$.
\end{lemma}

\begin{lemma}\label{lem:local commutativity version 2 for d}
Let $i\geq 1$. Then $\down_{i}\down_{i+1}^2=\down_{i+1}\down_{i}\down_{i+1}$.
\end{lemma}

\begin{lemma}\label{lem:commutativityjdtanddown}
Let $i\neq j$ be positive integers. Then
\begin{align*}
\addu_i\down_j=\down_j\addu_i.
\end{align*}
\end{lemma}
\begin{proof}
Let us first consider the case $1\leq i\leq j-1$. Then by Lemmas \ref{lem:commutativity of a and d distinct} and \ref{lem:non local commutativity for d}, we have that $\down_j$ commutes with $a_i,\down_1,\ldots,\down_{i-1}$. Hence $\addu_i\down_j=\down_j\addu_i$ in this case.

Now consider the case where $i>j\geq 1$. Then $\down_j\addu_i=\down_ja_i\down_1\down_2\cdots \down_{i-1}$. Again, using Lemmas \ref{lem:commutativity of a and d distinct} and \ref{lem:non local commutativity for d}, we can write this as 
\begin{align*}a_i\down_1\cdots \down_{j-2}\down_j\down_{j-1}\down_j\cdots \down_{i-1}.\end{align*}
Using Lemma \ref{lem:local commutativity version 2 for d}, we can write the above as 
\begin{align*}
a_i\down_1\cdots \down_{j-2}\down_{j-1}\down_{j}\down_j\down_{j+1}\cdots \down_{i-1}.
\end{align*}
Notice at this stage, if we assume $j=i-1$, then we have shown that $\addu_i\down_j=\down_j\addu_i$. So let us assume $i-j\geq 2$. Using Lemma \ref{lem:local commutativity version 1 for d}, we can transform the above expression to 
\begin{align*}
a_i\down_1\cdots \down_{j-2}\down_{j-1}\down_{j}\down_{j+1}\down_{j}\down_{j+2}\cdots \down_{i-1}.
\end{align*}
Now  Lemma \ref{lem:non local commutativity for d} easily establishes that the above expression equals
\begin{align*}
a_i\down_1\cdots \down_{j-2}\down_{j-1}\down_{j}\down_{j+1}\down_{j+2}\cdots \down_{i-1}\down_j
\end{align*}
and we are done.
\end{proof}

\begin{lemma}\label{lem:equalcommutativity}
Let $i\geq 1$. Then $\addu_i\down_i=\down_{i+1}\addu_{i+1}$.
\end{lemma}

\begin{proof}
Notice that $\addu_i\down_i=a_i\downrow_{[i]}$. Furthermore, Lemma \ref{lem:d and a same} states that $a_i=\down_{i+1}a_{i+1}$, and hence $\addu_i\down_i=\down_{i+1}a_{i+1}\downrow_{[i]}$. Since $\addu_{i+1}=a_{i+1}\downrow_{[i]}$, by definition, the claim follows.
\end{proof}

%%%relations t and d

\subsection{Relations satisfied by operators of type $\addt$ and $\down$}\label{subsec:td} We now give two useful lemmas, but this time regarding the box adding and box removing operators. Again, if desired, this subsection can be safely skipped for now and referred to later. In all the proofs we assume that $\alpha$ is a weak composition.

\begin{lemma}\label{lem: commutativity of t_i and d_j}
Let $i\neq j$ be positive integers. Then \begin{align*}
\addt_{i}\down_{j}=\down_{j}\addt_i.
\end{align*}
\end{lemma}
\begin{proof}
Let $\alpha = (\alpha_1,\ldots,\alpha_{\ell})$. 
First consider the case $i=1$. If $\alpha$ does not have a part equalling $j$, then $\addt_1\down_j(\alpha)=0$. Note now that, since $j\neq 1$, we have that $\down_j\addt_1(\alpha)=\down_j((1,\alpha_1,\ldots, \alpha_{\ell}))=0$ as well.

Hence we can assume that $i\geq 2$. If $\alpha$ does not have a part equalling $i-1$, then using the fact that $i\neq j$, we get that $\down_j(\alpha)$ does not have a part equalling $i-1$ either (assuming it does not equal $0$ already). This implies that $\addt_i\down_j(\alpha)=0$. Our assumption that $\alpha$ has no part equalling $i-1$ also implies that $\down_j\addt_i(\alpha)=0$.

Finally assume that $\alpha$ does have a part equalling $i-1$, and let $\alpha_r$ denote the leftmost such part. Then $$\addt_i(\alpha)=(\alpha_1,\ldots, \alpha_{r}+1,\ldots,\alpha_{\ell}).$$ If $\alpha$ does not have a part equalling $j$, then neither does $\addt_i(\alpha)$. This follows from the fact that $i\neq j$. This immediately implies that $\addt_i\down_j(\alpha)=\down_j\addt_i(\alpha)=0$ in this case. If $\alpha$ does have a part equalling $j$, then let $\alpha_s$ denote the rightmost such part. Note that $\alpha_s$ continues to be the rightmost part equalling $j$ in $\addt_i(\alpha)$ as well (unless there is a single part equalling $j=i-1$, in which case $\addt _i\down _j (\alpha) = \down _j \addt _i (\alpha) = 0$). Again, this follows since $i\neq j$. Thus we get that $$\addt_i\down_j(\alpha)=\down_j\addt_i(\alpha)=(\alpha_1,\ldots,\alpha_r+1,\ldots, \alpha_s-1,\ldots,\alpha_{\ell})$$if $r<s$ and
$$\addt_i\down_j(\alpha)=\down_j\addt_i(\alpha)=(\alpha_1,\ldots,\alpha_s-1,\ldots, \alpha_r+1,\ldots,\alpha_{\ell})$$if $s<r$.
\end{proof}

The proof of the next lemma consists of a number of small case analyses that are similar in style to the above proof, and hence we omit them for brevity.

\begin{lemma}\label{lem:zero contribution}
Let $i$ be a positive integer. When $i=1$ we have the following.
\begin{enumerate}
\item If $\alpha$ has parts equalling 1, then $\down_1\addt_1(\alpha)=\addt_1\down_1(\alpha) \neq 0.$
\item If $\alpha$ has no parts equalling 1, then $\down_1\addt_1(\alpha)=\alpha$ and $\addt_1\down_1(\alpha) = 0.$
\end{enumerate}
When $i\geq 2$ we have the following.
\begin{enumerate}
\item If $\alpha$ has parts equalling both $i$ and $i-1$, then $\down_i\addt_i(\alpha)=\addt_i\down_i(\alpha) \neq 0.$
\item If $\alpha$ has parts equalling $i$ and no parts equalling $i-1$, then $\down_i\addt_i(\alpha)=0$ and $\addt_i\down_i(\alpha) =\alpha.$
\item If $\alpha$ has no parts equalling $i$ and parts equalling $i-1$, then $\down_i\addt_i(\alpha)=\alpha$ and $\addt_i\down_i(\alpha) =0.$
\item If $\alpha$ has parts equalling neither $i$ nor $i-1$, then $\down_i\addt_i(\alpha)=\addt_i\down_i(\alpha) = 0.$
\end{enumerate}
In particular, if $\down_i\addt_i(\alpha)$ and $\addt_i\down_i(\alpha)$ are nonzero, then $\down_i\addt_i(\alpha)=\addt_i\down_i(\alpha)$.
\end{lemma}

%%% Dual graphs
\section{Dual graphs from composition posets}\label{sec:dualgraphs} 
We now recall terminology pertaining to graded graphs and filtered graphs, and follow the notation of \cite{patrias-pyl}.
Let $G$ be a graph consisting of a set of vertices $P$ endowed with a rank function $\rho:P\rightarrow \mathbb{Z}$ with vertices $x,y\in P$ and $y$ is of rank weakly greater than $x$. Then $G$ is called a \emph{graded graph} when the edge set $E$ satisfies if $(x,y)\in E$ then $\rho(y)= \rho(x) + 1$. The graph $G$ is called a \emph{weak filtered graph} when the edge set $E$ satisfies if $(x,y)\in E$ then $\rho(y)\geq \rho(x)$, and a \emph{strong filtered graph} when the edge set $E$ satisfies if  $(x,y)\in E$ then $\rho(y)> \rho(x)$.  Now given a field $K$ of characteristic $0$, the vector space $KP$ is the space consisting of all formal linear combinations of vertices of $G$. Then we define the up and down operators $U,D \in End(KP)$ associated with $G$ to be
$$U(x) = \sum _y m(x,y) y$$
$$D(y) = \sum _x m(x,y) x$$where $x$ and $y$ are vertices of $G$, $y$ is of weakly greater rank than $x$, and $m(x,y)$ is the number of edges connecting $x$ and $y$. With this in mind, let $G_1$ be a graded graph with up operator $U$ and $G_2$ be a graded graph with down operator $D$ such that $G_1$ and $G_2$ have a common vertex set $P$ and rank function $\rho$. Then $G_1$ and $G_2$ are \emph{dual graded graphs} if and only if on $KP$
\begin{equation}\label{eq:du-ud=i}
DU-UD=\Id
\end{equation}where $\Id$ is the identity operator on $KP$. Similarly let $\widetilde{G_1}$ be a weak filtered graph with up operator $\widetilde{U}$ and $\widetilde{G_2}$ be a strong filtered graph with down operator $\wD$ such that $\widetilde{G_1}$ and $\widetilde{G_2}$ have a common vertex set $P$ and rank function $\rho$. Then $\widetilde{G_1}$ and $\widetilde{G_2}$ are \emph{dual filtered graphs} if and only if on $KP$
\begin{equation}\label{eq:wdu-uwd=wd+i}
\wD \widetilde{U}-\widetilde{U}\wD =\wD + \Id .
\end{equation}

%%%Rc

\subsection{Dual graphs and the right composition poset}\label{subsec:DRc} Observe that our composition posets $\Rc$ and $\Qc$ defined in Subsection~\ref{subsec:posets} with vertex set being the set of all compositions, whose rank function is given by the size of a composition and whose  edge sets are the respective cover relations, are both clear examples of graded graphs. By the definition of the cover relation $\coverr$ it follows that the up operator associated with $\Rc$ is given by
\begin{equation}\label{eq:URc}
U= \sum _{i\geq 1} \addu _i.
\end{equation}

\begin{example}\label{ex:URc}
Let $\alpha $ be the composition $(2,1,3)$. Then 
\begin{align*}U((2,1,3))&=(2,1,3,1)+(2,0,3,2)+(1,0,3,3)+(2,1,0,4)\\
&=(2,1,3,1)+(2,3,2)+(1,3,3)+(2,1,4).
\end{align*}
\end{example}
Similarly, by the definition of the cover relation $\coverq$ it follows that the down operator associated with $\Qc$ is given by
\begin{equation}\label{eq:DQc}
D= \sum _{i\geq 1} \down _i.
\end{equation}

\begin{example}\label{ex:DQc}
Let $\alpha $ be the composition $(2,1,3)$. Then by Example~\ref{ex:down}
$$D((2,1,3))= (2,0,3)+(1,1,3)+(2,1,2)= (2,3)+(1,1,3)+(2,1,2).$$
\end{example}

Moreover we have the following.

\begin{theorem}\label{the:RcQc}
$\Rc$ and $\Qc$ are dual graded graphs, that is, on compositions
$$DU-UD=\Id.$$
\end{theorem}

\begin{proof}
Notice that 
\begin{align*}
DU=\sum_{\substack{i\neq j\\i,j\geq 1}}\down_j\addu_i+\sum_{k\geq 1}\down_k\addu_k
\end{align*}
and 
\begin{align*}
UD=\sum_{\substack{i\neq j\\i,j\geq 1}}\addu_i\down_j+\sum_{k\geq 1}\addu_k\down_k.
\end{align*}
Using Lemma \ref{lem:commutativityjdtanddown} and Lemma \ref{lem:equalcommutativity}, we reach the conclusion that 
\begin{align*}
DU-UD=\down_1\addu_1.
\end{align*}
By Lemma \ref{lem:d and a same}, $\down_1\addu_1=a_0=Id$. This finishes the proof.
\end{proof}

\begin{example}\label{ex:RcQc}
Let $\alpha=(2,1,3)$. Then suppressing commas and parentheses for ease of comprehension, we have by Examples~\ref{ex:URc} and \ref{ex:DQc} that
\begin{align*}
DU(\alpha)&=D(2131+2032+1033+2104)\\&=2130+1131+2121+2031+2022%\\&&
+0033+1032+2004+1104+2103
%\\&=& 213+1131+2121+231+222\\&&+33+132+24+114+213.
\end{align*}
and
\begin{align*}
UD(\alpha)&=U(203+113+212)\\&=2031+0033+2004+1131+1032%\\&
+1104+2121+2022+2103.
%\\&=& 231+33+24+1131+132\\&&+114+2121+222+213.
\end{align*}
Thus $(DU-UD)(\alpha)=213=\Id (\alpha)$. 
\end{example}

To describe our results in the context of dual filtered graphs, we need the following.
\begin{definition}\label{def:Qct} 
Let $\Qct$ be the graded graph whose vertex set is the set of all compositions, whose rank function is given by the size of a composition, and whose edge set is as follows. Given compositions $\alpha$ and $\beta$ such that the size of $\alpha$ is strictly greater than $\beta$, we have the edge 
$$(\beta,\alpha) \mbox{ if and only if } \down _I(\alpha) = \beta$$for some finite $\emptyset \neq I \subset \bN$. 
\end{definition}
As before, when computing $\down_I(\alpha)$ in Definition~\ref{def:Qct}, we ignore all parts that equal $0$.

\begin{example}\label{ex:Qct}We have an edge between  $\beta = (4,1,3,1,1)$ and $\alpha = (4,1,4,2,1)$ in $\Qct$ since $\down _{\{2,4\}} (\alpha) = \beta$.
\end{example}

\begin{remark}
Observe that the relation $\coverqt$ on compositions defined by $\beta\coverqt \alpha$ if and only if $\beta=\down_I(\alpha)$ does not give rise to a poset structure, since transitivity is not satisfied. For example,  $\down _{\{1,4\}} ((4,1,4,1))=(4,1,3)$ and $\down _{\{1,4\}} ((4,1,3))=(3,3)$, but no $I$ exists such that $\down _I ((4,1,4,1))=(3,3)$.%This is in contrast to the other graded graphs constructed in this article.
\end{remark}

Clearly, we have that $\Qct$  is an example of a strong filtered graph by definition. The associated down operator is given by
\begin{equation}\label{eq:RcQct}
\wD = \sum _{I\subset \bN} \downrow _I
\end{equation}where the sum is over all finite but nonempty subsets of $\bN$. Hence we can relate $\Rc$ and $\Qct$ as follows, since any graded graph, such as $\Rc$, is also a weak filtered graph.

\begin{theorem}\label{the:RcQct}
$\Rc$ and $\Qct$ are dual filtered graphs, that is, on compositions
$$\wD U-U\wD =\wD + \Id.$$
\end{theorem}

\begin{proof}
First note that the operator $\wD U$ has the following expansion.
\begin{align*}
\wD U&=\sum_{\substack{i\geq 1\\I\subset \bN}}\downrow_I\addu_i\\&=\sum_{\substack{I\subset \bN\\i\in I}}\downrow_I\addu_i + \sum_{\substack{I\subset \bN\\i\geq 1,i\notin I}}\downrow_I\addu_i
\end{align*} 
In a similar manner, we obtain the following expansion for $U\wD $.
\begin{align*}
U\wD &=\sum_{\substack{i\geq 1\\I\subset \bN}}\addu_i\downrow_I\\&=\sum_{\substack{I\subset \bN\\i\in I}}\addu_i\downrow_I + \sum_{\substack{I\subset \bN\\i\geq 1,i\notin I}}\addu_i\downrow_I
\end{align*}
Using Lemma \ref{lem:commutativityjdtanddown}, we obtain that
\begin{align*}
\wD U-U\wD=\sum_{\substack{I\subset \bN\\i\in I}}\downrow_I\addu_i-\sum_{\substack{I\subset \bN\\i\in I}}\addu_i\downrow_I.
\end{align*}

Now consider a fixed set $I\subset \bN$ and $i\in I$. We will next show that the operator $\downrow_I\addu_i$ corresponds to either to a unique operator $\addu_{i'}\downrow_{I'}$ where $i' \in I'$, or an operator $a_0\downrow_{I'}$ where $I'$ might be the empty set. 

Let $j\in I$ be the smallest positive integer such that $j-1\notin I$ but every integer $k$ satisfying $j\leq k\leq i$ belongs to $I$. Consider the following sets.
\begin{align*}
A&=\{k\suchthat k\in I, k <j\}\\
B&=\{k \suchthat j\leq k\leq i\}\\
C&=\{k \suchthat k\in I, k>i\}
\end{align*}
Clearly, we have that $I=A\amalg B\amalg C$ where $\amalg$ denotes disjoint union. Define the set $I'$ to be $A\amalg (B-1) \amalg C$. Notice that $I'$ can be the empty set (precisely in the case where $A$ and $C$ are empty, while $B=\{1\}$). Now we have the following sequence of equalities using Lemma \ref{lem:commutativityjdtanddown} and Lemma \ref{lem:equalcommutativity}.
\begin{align*}
\downrow_I\addu_i&=\downrow_{A}\downrow_{B}\downrow_{C}\addu_i\\&=\downrow_A\downrow_B\addu_i\downrow_C\\&=\downrow_A\addu_{j-1}\downrow_{B-1}\downrow_C\\&=\addu_{j-1}\downrow_A\downrow_{B-1}\downrow_C\\&=\addu_{j-1}\downrow_{I'}
\end{align*}
Given the invertibility of our computation,  it is clear how to recover $\downrow_I\addu_i$ starting from $\addu_{j-1}\downrow_I'$. Furthermore, if $j\neq 1$, then we clearly have that $j-1\in I'$. The above thus implies that
\begin{align*}
\sum_{\substack{I\subset \bN\\i\in I}}\downrow_I\addu_i-\sum_{\substack{I\subset \bN\\i\in I}}\addu_i\downrow_I &= a_0+a_0\wD
\end{align*}
thereby finishing the proof.
\end{proof}

\begin{example}\label{ex:RcQct}
Let $\alpha = (1,2)$. Then suppressing commas and parentheses as before, we have that
$$\wD (\alpha)= (02+11+10).$$
Therefore
\begin{align*}
\wD U (\alpha) &= \wD (121+022+103)\\
&=120+111+110+021+020%\\&&
+003+102+002+101+100
\end{align*}and
\begin{align*}
U \wD  (\alpha) &= U (02+11+10)\\
&= 021+003+111+102+101+002.\end{align*}
Thus $(\wD U-U\wD)(\alpha)=2+11+1+12=(\wD + \Id) (\alpha)$. 
\end{example}

\begin{remark}\label{rem:Fomin}
It is worth noting the connection between our results here and Fomin's work in \cite{fomin-schur}. In particular, note that the relations \cite[Equation 1.9]{fomin-schur} satisfied by his box adding and box removing operators on partitions (denoted therein by $u$ and $d$, respectively) are the same as those satisfied by the jdt operators and box removing operators on compositions. The relations are easy to establish in the case of partitions, but as we have seen,  deriving the same relations in the case of compositions is  more delicate. 

Fomin then uses these operators to define generating functions $A(x)$ and $B(y)$ that add or remove horizontal strips in all possible ways respectively, and then uses \cite[Equation 1.9]{fomin-schur} to prove the following commutation relation \cite[Theorem 1.2]{fomin-schur}.
\begin{align*}
A(x)B(y)=B(y)A(x)(1-xy)^{-1}
\end{align*} 
He later notes that the dual graded graph nature of Young's lattice is encoded in the aforementioned identity. More precisely it follows from comparing the coefficient of $xy$ on either side \cite[Equation 1.13]{fomin-schur}. In fact, one can obtain various identities by comparing coefficients and can verify that the relations describing dual filtered graphs can be obtained by setting $y=1$ and then subsequently comparing the coefficient of $x$ on either side.
Thus in a sense, the relations uniformly establish both the dual graded graph and the dual filtered graph structures on Young's lattice and $\Rc$.
\end{remark}

We now proceed to discuss $\Lc$ defined using box adding operators. We will establish that this poset can also be endowed with a structure of a dual graded graph and a dual filtered graph. But the relations satisfied in this case are different than the ones we have encountered, and we cannot use Fomin's commutation relation in this setting. In fact, as we will see, the cancellations in the case of $\Lc$ are more subtle despite the  simplicity of the action of $\addt$ compared to the action of $\addu$.

%%%Lc
\subsection{Dual graphs and the left composition poset}\label{subsec:DLc} Our composition poset $\Lc$ with vertex set being the set of all compositions,  whose rank function is given by the size of a composition and whose edge set is the cover relations is clearly a graded graph and hence also a weak filtered graph. By the definition of the cover relation $\coverc$ it follows that the up operator associated with $\Lc$ is given by
\begin{equation}\label{eq:ULc}
\Ut = \sum _{i\geq 1} \addt _i.
\end{equation}

\begin{example}\label{ex:ULc}
Let $\alpha $ be the composition $(2,1,3)$. Then
$$\Ut((2,1,3))=(1,2,1,3)+(2,2,3)+(3,1,3)+(2,1,4).
$$
\end{example}

Again $\Qc$ and $\Qct$ are respectively a graded graph and a strong filtered graph with respective down operators $D$ and $\wD$.

%%%setup
For the remainder of this section, we will fix a composition $\alpha$.
This given, we will define a function $\Phi: Y\rightarrow X$, where the sets $X$ and $Y$ are defined as follows.
\begin{align*}
X&=\{\downrow_I\addt_i\suchthat I\subset \bN, i\in I, \downrow_I\addt_i(\alpha)\neq 0 \}\\
Y&=\{\addt_i\downrow_I\suchthat I\subset \bN, i\in I, \addt_i\downrow_I(\alpha)\neq 0\}
\end{align*} 
Consider $w=\addt_i\downrow_I\in Y$. 
Decompose $I=A\amalg\{i\}\amalg B$ where
\begin{align*}
A&=\{j\in I\suchthat j<i\}\\
B&=\{j\in I\suchthat j>i\}.
\end{align*}
By Lemma \ref{lem: commutativity of t_i and d_j}, we have that $w=\downrow_A\addt_i\down_i\downrow_B$.
Let $k$ denote the largest part of $\alpha$ that is strictly less than $i$. Then $k\geq \max(A)$ as follows. 

Decompose $A = A' \coprod \{m\}$ with $m = \max(A)$. Then we have  $\mathfrak d_{A'} \mathfrak d_{m} \mathfrak t_i \mathfrak d_i \mathfrak d_B (\alpha) = w(\alpha) \neq 0$ and it follows that $m$ is a part in the composition $\mathfrak t_i \mathfrak d_i \mathfrak d_B (\alpha)$ (otherwise, $\mathfrak d_{m}(\mathfrak t_i \mathfrak d_i \mathfrak d_B (\alpha)) = 0$, which implies that $w(\alpha) = 0$ contradicting $w(\alpha) \neq 0$). Hence the largest part of $\alpha$ strictly less than $i$ is at least $m = \max(A)$.

If such a part does not exist, we define $k$ to be $0$. Let $i'=k+1$. 
Now let $I'=A\amalg\{i'\}\amalg B$ and $$\Phi(w)=\downrow_{I'}\addt_{i'}=\downrow _A \down _{i'}\addt _{i'} \downrow _B.$$

Then the following can be proved using Lemmas~\ref{lem: commutativity of t_i and d_j} and  \ref{lem:zero contribution}.

\begin{lemma}\label{lem: properties of Phi}
Let $w=\addt_i\downrow_I=\addt_i\downrow_A\down_i\downrow_B = \downrow_A\addt_i\down_i\downrow_B\in Y$ and let $w'=\Phi(w)$. Then the following statements hold.
\begin{enumerate}
\item $w'(\alpha)=w(\alpha)$ if $i=1$.
\item  $w'(\alpha)=w(\alpha)$ if $i\geq 2$ and $i$ is not the smallest part of $\downrow_B(\alpha)$.
\item $w'(\alpha)=(0,w(\alpha))$ if $i\geq 2$ and $i$ is the smallest part of $\downrow_B(\alpha)$.
\end{enumerate}
In particular, $\Phi: Y\rightarrow X$ and, at the level of compositions, we have that   $\Phi(w)(\alpha)=w(\alpha)$ for all $w\in Y$.
\end{lemma}

The next step for us is to identify the image of $Y$ under the map $\Phi$. The image of $Y$ is a very special subset of $X$, which has the following explicit description. Let the largest part of $\alpha$ be $m$. Define $Z$ as follows.
\begin{align*}
Z=\{\downrow_I\addt_i\in X\suchthat i\leq m\}
\end{align*}
Thus in other words, $Z$ is the subset comprising of words that never add a box to the largest part. Note that by the definition of $\Phi$ we have that $\Phi(Y)\subseteq Z$ since if $w\in Y$ and $\Phi(w)$ has rightmost operator $\addt _j$ then $j\leq m$. Our next aim is to  find the inverse of $\Phi$.

Consider $w=\downrow_I\addt_i\in Z$. Writing $I=A\amalg\{i\}\amalg B$ in the usual way, and using Lemma \ref{lem: commutativity of t_i and d_j} allows us to write $w$ as shown below.
\begin{align*}
w=\downrow_A\down_i\addt_i\downrow_B
\end{align*}
Let $i''$ be the smallest part of $\downrow_B(\alpha)$ weakly greater than $i$. This always exists by our hypothesis that $w\in Z$. We define $\Psi (w)$ to be
\begin{align*}
\Psi(w)&=\downrow_A\addt_{i''}\down_{i''}\downrow_B\\
&=\addt_{i''}\downrow_{I''}
\end{align*}
where $I''=A\amalg \{i''\}\amalg B$. It is straightforward to see that if $k$ is the largest part of $\alpha$ strictly less than $i$, $i'=k+1$, and $i''$ is the smallest part of $\downrow _B(\alpha)$ weakly greater than $i'$, then $i''=i$ and hence
$$\Psi(\Phi(w))= \downrow_A\addt_i\down_i\downrow_B=w$$so $\Psi$ is the inverse of $\Phi$. Hence we have the bijection
\begin{equation}\label{eq:PhiofY}
\Phi(Y)=Z.
\end{equation}

\begin{example}\label{ex:setup}
Consider the composition $\alpha=(2,6,1,4)$ and let $w=\addt_4\downrow_{\{1,4,5,6\}}$. Then $w (\alpha) = (2,4,0,4)$ so $w\in Y$. We have the following decomposition for $w$.
\begin{align*}
w=\downrow_{\{1\}}\addt_4\down_4\downrow_{\{5,6\}}
\end{align*}
Then the corresponding $A$, $B$ and $i$ are $\{1\}$, $\{5,6\}$ and $4$ respectively. Our method for constructing $\Phi(w)$ requires that first we find the largest part $k$ strictly less than $i$ in $\alpha$.  So it follows that $k=2$. This implies that
\begin{align*}
\Phi(w)=\downrow_{\{1,3,5,6\}}\addt_3 = \downrow _{\{1\}}\down _3\addt_3 \downrow _{\{5,6\}}
\end{align*}
and hence $\Phi(w)(\alpha)=(2,4,0,4)=w(\alpha)$ and   $\Phi(w)\in Z$. Lastly note that since  $\downrow_B(\alpha)=(2,4,1,4)$ we have for $\Phi(w)$ that its $i''=4$ and
$$\Psi(\Phi(w))=\Psi(\downrow_{\{1\}}\down _3\addt _3 \downrow _{\{5,6\}})=\downrow_{\{1\}}\addt _4 \down _4 \downrow _{\{5,6\}}=w$$as desired.\end{example}

Since $\Psi$ is the inverse of $\Phi$ we also have the following.

\begin{corollary}\label{cor:injectYandX}
$\Phi$ is an injection from $Y$ to $X$.
\end{corollary}

Consider the sets $P$ and $Q$ defined as follows.
\begin{align*}
P=\{\downrow_i\addt_i\suchthat i\geq 1, \downrow_i\addt_i(\alpha)\neq 0 \}\\
Q=\{\addt_i\downrow_i\suchthat i\geq 1, \addt_i\downrow_i(\alpha)\neq 0\}\\
\end{align*} 
Clearly, $P\subset X$ and $Q\subset Y$. Furthermore, we have that $\Phi(Q)$ maps into $P$. In fact, a stronger claim holds from the discussion prior to this: $$\Phi(Q)=P\setminus \{\down_{m+1}\addt_{m+1}\}$$where $m$ is the largest part of $\alpha$.

%%%

Then utilising all of the above we have the following two theorems.

\begin{theorem}\label{the:LcQc}
$\Lc$ and $\Qc$ are dual graded graphs, that is, on compositions
$$D\Ut -\Ut D=\Id.$$
\end{theorem} 

\begin{proof}
Firstly note that $D\Ut$ corresponds to the following expansion.
\begin{align*}
D\Ut&=\sum_{i,j\geq 1}\down_i\addt_j=\sum_{i,j\geq 1, i\neq j}\down_i\addt_j + \sum_{k\geq 1}\down_k\addt_k
\end{align*}
Also the operator $\Ut D$ corresponds to the expansion below.
\begin{align*}
\Ut D&=\sum_{i,j\geq 1}\addt_j\down_i=\sum_{i,j\geq 1, i\neq j}\addt_j\down_i + \sum_{k\geq 1}\addt_k\down_k
\end{align*}
Then, on using Lemma \ref{lem: commutativity of t_i and d_j}, we obtain the following.
\begin{align}\label{eq:D\Ut-\Ut D lem:d and a same}
D\Ut -\Ut D&=\sum_{k\geq 1}\down_k\addt_k-\sum_{k\geq 1}\addt_k\down_k
\end{align}
Taking $\alpha$ into account we can rewrite the above equation as  stating the following.
\begin{align*}
(D\Ut-\Ut D)(\alpha)=\sum_{w\in P}w(\alpha)-\sum_{w\in Q}w(\alpha)=\down_{m+1}\addt_{m+1}(\alpha)+\sum_{w\in Q}(\Phi(w)-w)(\alpha)
\end{align*}
Now  at the level of compositions we have $\sum_{w\in Q}(\Phi(w)-w)(\alpha)=0$ by Lemma \ref{lem: properties of Phi}, and $\down_{m+1}\addt_{m+1}(\alpha)=\alpha$. This implies the claim.
\end{proof}

\begin{example}\label{ex:LcQc}
Let $\alpha=(2,1,3)$. Then suppressing commas and parentheses as before, we have that
\begin{align*}
D\Ut(\alpha)&=D(1213+223+313+214)\\&=1203+1113+1212+213+222%\\&&
+303+312+204+114+213
\end{align*}
and
\begin{align*}
\Ut D(\alpha)&= \Ut(203+113+212)\\&= 1203+303+204+1113+213%\\&&
+114+1212+222+312.
\end{align*}
Thus $(D\Ut -\Ut D)(\alpha)=213=\Id (\alpha)$. 
\end{example}

\begin{theorem}\label{the:LcQct}
$\Lc$ and $\Qct$ are dual filtered graphs, that is, on compositions
$$\wD \Ut-\Ut\wD =\wD + \Id.$$
\end{theorem}

\begin{proof}
The beginning of the proof is very similar to that in Theorem \ref{the:RcQct} but with $\addt_i$ instead of $\addu _i$. Using Lemma \ref{lem: commutativity of t_i and d_j} we obtain the following equality.
\begin{align*}
\wD \Ut-\Ut\wD=\sum_{\substack{I\subset \bN\\i\in I}}\downrow_I\addt_i-\sum_{\substack{I\subset \bN\\ i\in I}}\addt_i\downrow_I
\end{align*}

Now for the fixed composition $\alpha$, we can rewrite the above equation as follows.
\begin{align}\label{eq: equation before cancellation}
(\wD\Ut-\Ut\wD)(\alpha)&=\sum_{w\in X}w(\alpha)-\sum_{w\in Y}w(\alpha)\nonumber\\
&=
\sum_{w\in X\setminus Z}w(\alpha)+\sum_{w\in Z}w(\alpha)-\sum_{w\in Y}w(\alpha)
\end{align}
At the level of  compositions, Lemma \ref{lem: properties of Phi} implies that 
\begin{align*}
\sum_{w\in Y}(\Phi(w)(\alpha)-w(\alpha))=0.
\end{align*}
Using the above and Equation \eqref{eq:PhiofY} in Equation \eqref{eq: equation before cancellation}   at the level of   compositions gives
\begin{align*}
(\wD \Ut -\Ut\wD)(\alpha)=\sum_{w\in X\setminus Z}w(\alpha).
\end{align*}
Observe now that every element of $X\setminus Z$ has the form $\downrow_A\down_{m+1}\addt_{m+1}$ where $A$ consists only of instances of $\down_i$ where $i\leq m$ and $m$ is the largest part of $\alpha$. Furthermore we do have the possibility that $A$ is empty. Additionally, it is easy to see that $\down_{m+1}\addt_{m+1}$ is the identity map. The preceding discussion allows us to conclude the following equality at the level of compositions, thereby finishing the proof.
\begin{align*}
(\wD\Ut -\Ut\wD)(\alpha)=(\wD+\Id)(\alpha)
\end{align*}\end{proof}

\begin{example}\label{ex:LcQct}
Let $\alpha = (1,2)$. Then suppressing commas and parentheses as before, we have that
$$\wD (\alpha)= (02+11+10).$$
Therefore
\begin{align*}
\wD \Ut (\alpha) &= \wD (112+22+13)\\
&=102+111+110+21+20%\\&&
+ 03+12+02+11+10
\end{align*}and
\begin{align*}
\Ut \wD  (\alpha) &= \Ut (02+11+10)\\
&= 102+03+111+21+110+20.\end{align*}
Thus $(\wD \Ut-\Ut\wD)(\alpha)=2+11+1+12=(\wD + \Id) (\alpha)$. 
\end{example}

\section*{Acknowledgements} The author would like to thank Vasu Tewari for many enjoyable and enlightening conversations. She is also most grateful to the referee for their deep and thoughtful feedback, which she implemented if journal space permitted and used for editing guidance if not.
%: Bibliography
%\bibliography{SQScoV3}

\begin{thebibliography}{99}

\bibitem{berg-saliola-serrano}
{\sc C.~Berg, F.~Saliola and L.~Serrano},
{\em The down operator and expansions of near rectangular $k$-Schur functions},
J. Combin. Theory Ser. A 120 (2013) 623--636.

\bibitem{bergeron-lam-li}
{\sc N.~Bergeron, T.~Lam and H.~Li},
{\em{Combinatorial Hopf algebras and towers of algebras -- dimension, quantization and functorality}},
 Algebr. Represent. Theory 15 (2012) 675--696.
 
 
\bibitem{BLvW} 
{\sc C.~Bessenrodt, K.~Luoto and S.~van Willigenburg}, 
{\em Skew quasisymmetric Schur functions and noncommutative Schur functions},
Adv. Math. 226 (2011) 4492--4532.


   
   \bibitem{DKLT} 
{\sc G.~Duchamp, D.~Krob, B.~Leclerc and J-Y.~Thibon}, 
{\em Fonctions quasi-sym\'etriques, 
fonctions sym\'etriques non-commutatives, et alg\`ebres de Hecke \`a $q = 0$},
C. R. Math. Acad. Sci. Paris
 322 (1996) 
  107--112. 
 
  
\bibitem{fomin-dual1}
{\sc S.~Fomin},
{\em Duality of graded graphs},
J. Algebraic Combin. 3 (1994) 357--404.

\bibitem{fomin-dual2}
{\sc S.~Fomin},
{\em Schensted algorithms for dual graded graphs},
J. Algebraic Combin. 4 (1995) 5--45.

  \bibitem{fomin-schur}
{\sc S.~Fomin},
{\em {Schur operators and Knuth correspondences}},
J. Combin. Theory Ser. A 72 (1995) 277--292.

\bibitem{gaetz}
{\sc C.~Gaetz},
{\em {Dual graded graphs and Bratelli diagrams of towers of groups}},
\url{arXiv:1803.11168}  


\bibitem{GKLLRT}
{\sc I.~Gelfand, D.~Krob, A.~Lascoux, B.~Leclerc, V.~Retakh and J.-Y. Thibon},
  {\em {Noncommutative symmetric functions}}, 
  Adv. Math. 112 (1995) 218--348.
  

\bibitem{gessel}
{\sc I.~Gessel}, {\em {Multipartite P-partitions and inner products of skew
  Schur functions}},
 {Combinatorics and algebra, Proc. Conf., Boulder/CO 1983, Contemp.
  Math. 34 (1984) 289--301}.


\bibitem{QS}
{\sc J.~Haglund, K.~Luoto, S.~Mason and S.~van Willigenburg}, 
{\em {Quasisymmetric {S}chur functions}}, 
J. Combin. Theory Ser. A 118 (2011) 463--490. 

  \bibitem{hersh-hsiao} 
{\sc P.~Hersh and S.~Hsiao},
{\em Random walks on quasisymmetric functions},
Adv. Math.
222 (2009)
 782--808.

  
\bibitem{lam-sign} 
 {\sc T.~Lam},
 {\em Signed differential posets and sign-imbalance},
J. Combin. Theory Ser. A 115 (2008) 466--484.
 
 \bibitem{lam}
 {\sc T.~Lam},
 {\em Quantized dual graded graphs},
 Electron. J. Combin. 17 (2010).
 
 \bibitem{lam-shimozono}
 {\sc T.~Lam and M.~Shimozono},
 {\em Dual graded graphs for Kac-Moody algebras},
Algebra Number Theory 1 (2007) 451--488.

\bibitem{miller}
{\sc A.~Miller}, {\em {Differential posets have strict rank growth: a conjecture of Stanley}}, 
Order 30 (2013) 657--662.

\bibitem{nuts}
{\sc J.~Nzeutchap}, {\em Dual graded graphs and Fomin's $r$-correspondences associated to the Hopf
algebras of planar binary trees, quasi-symmetric functions and
noncommutative symmetric functions}, FPSAC 2006.


  \bibitem{patrias-pyl}
{\sc R.~Patrias and P.~Pylyavskyy},
{\em {Dual filtered graphs}},
\url{arXiv:1410.7683}  


 \bibitem{stanley-diff}
{\sc R.~Stanley}, {\em Differential posets}, J. Amer. Math. Soc. 1 
(1988) 919--961.

\bibitem{stanley-var}
{\sc R.~Stanley}, {\em Variations on differential posets}, IMA Vol. Math. Appl. 19 (1990) 145--165.
  
  
 \bibitem{stanley-zanello}
 {\sc R.~Stanley and F.~Zanello}, \emph{On the rank function of a differential poset}, Electron. J. Combin. 19 (2012). 
   
  
\bibitem{tewari}
{\sc V.~Tewari},
{\em Backward jeu de taquin slides for composition tableaux and a noncommutative Pieri rule},
Electron. J. Combin. 22 (2015).

\end{thebibliography}
\bibliographystyle{amsplain}
%    Insert the bibliography data here.
%\def\cprime{$'$}

\end{document}